\documentclass[11pt]{amsart}
\usepackage{amssymb,latexsym,epsfig,color}
 
\newtheorem{theorem}{Theorem}[section]
\newtheorem{proposition}[theorem]{Proposition}
\newtheorem{corollary}[theorem]{Corollary}

\newtheorem{lemma}[theorem]{Lemma}

\newtheorem{problem}[theorem]{Problem}

\theoremstyle{definition}
\newtheorem{definition}[theorem]{Definition}
\newtheorem{remark}[theorem]{Remark}
\newtheorem{example}[theorem]{Example}
\newtheorem{algorithm}[theorem]{Algorithm}

\newtheorem{observation}[theorem]{Observation}

\def\NN{\ensuremath{\mathbb{N}}}
\def\ZZ{\ensuremath{\mathbb{Z}}}

\def\RR{\ensuremath{\mathbb{R}}}

\newcommand{\B}{{\mathcal B}}

\def\A{\ensuremath{\mathcal{A}}}

\def\U{\ensuremath{\mathcal{U}}}

\def\aaa{\ensuremath{{\bf{a}}}}
\def\bb{\ensuremath{{\bf{b}}}}

\def\dd{\ensuremath{{\bf{d}}}}
\def\ee{\ensuremath{{\bf{e}}}}

\def\hh{\ensuremath{{\bf{h}}}}
\def\kk{\ensuremath{{\bf{k}}}}
\def\mm{\ensuremath{{\bf{m}}}}

\def\rr{\ensuremath{{\bf{r}}}}
\def\sss{\ensuremath{{\bf{s}}}}
\def\ttt{\ensuremath{{\bf{t}}}}
\def\uu{\ensuremath{{\bf{u}}}}
\def\vv{\ensuremath{{\bf{v}}}}

\def\xx{\ensuremath{{\bf{x}}}}
\def\yy{\ensuremath{{\bf{y}}}}

\def\Hilb{\textup{Hilb}}

\def\max{\ensuremath{\textup{max}}}

\def\cone{\ensuremath{\textup{cone}}}
\def\conv{\ensuremath{\textup{conv}}}

\def\scr{\ensuremath{\textup{SCR}}}
\def\Frac{\ensuremath{\textup{FRAC}}}

\def\stab{\ensuremath{\textup{STAB}}}

\newcommand{\bigo}{{\mathcal O}}

\def\endproof{\hfill$\square$\medskip}

\title[Small Chv{\'a}tal Rank]{Small Chv{\'a}tal Rank} 
\author{Tristram Bogart}
\address{Dept. of Mathematics and Statistics, Queen's Univ, Kingston, ON K7L 3N6}
\email{bogart@mast.queensu.ca}
\author{Annie Raymond}
\address{Berlin Mathematical School, Technical University, Berlin, 10623}
\email{raymond@math.tu-berlin.de}
\author{Rekha Thomas}
\address{Department of Mathematics, Univ. of Washington, 
       Seattle, WA 98195-4350} 
\email{thomas@math.washington.edu}
\thanks{All authors were partially supported by NSF grant DMS-0401047 and
  the Robert R. and Elaine K. Phelps Endowment at the University of Washington} 
\date{\today}

\begin{document}

\begin{abstract}
We propose a variant of the Chv{\'a}tal-Gomory procedure that will
  produce a sufficient set of facet normals for the integer hulls of
  all polyhedra $\{ \xx \,:\, A \xx \leq \bb\}$ as $\bb$ varies. The
  number of steps needed is called the small Chv{\'a}tal rank (SCR) of
  $A$.  We characterize matrices for which SCR is zero via the notion
  of supernormality which generalizes unimodularity. SCR is studied in
  the context of the stable set problem in a graph, and we show that
  many of the well-known facet normals of the stable set polytope
  appear in at most two rounds of our procedure. Our results reveal a
  uniform hypercyclic structure behind the normals of many complicated
  facet inequalities in the literature for the stable set
  polytope. Lower bounds for SCR are derived both in general and for
  polytopes in the unit cube.
\end{abstract}
\maketitle

\section{Introduction} \label{sec:intro}
The study of integer hulls of rational polyhedra is a fundamental area
of research in integer programming. For a matrix $A \in \ZZ^{m \times
  n}$ and a vector $\bb \in \ZZ^m$, consider the polyhedron 
$$ Q_\bb := \{ {\bf x} \in \RR^n \,:\, A{\bf x} \leq \bb \} $$ 
and its integer hull 
$$ Q_\bb^I := \textup{convex hull}(Q_\bb \cap \ZZ^n).$$ 
The {\em Chv{\'a}tal-Gomory} procedure is an algorithm for computing $Q_\bb^I$ from $Q_\bb$. This method involves iteratively adding rounds of {\em cutting planes} to $Q_\bb$ until $Q_\bb^I$ is obtained. The {\em Chv{\'a}tal rank} of $A \xx \leq \bb$ is the minimum number of rounds of cuts needed in the
Chv{\'a}tal-Gomory procedure to obtain $Q_\bb^I$, and the {\em
  Chv{\'a}tal rank} of $A$ is the maximum of the Chv{\'a}tal ranks of
$A \xx \leq \bb$ as $\bb$ varies in $\ZZ^m$.

In this paper we fix a matrix $A \in \ZZ^{m \times n}$ of rank $n$ and
look at the more basic problem of finding just the normals of a
sufficient set of inequalities that will cut out all integer hulls
$Q_\bb^I$ as $\bb$ varies in $\ZZ^m$. Given $A$, it is known that
there exists a matrix $M$ such that for each $\bb$, $Q_\bb^I = \{ \xx
\in \RR^n \,:\, M \xx \leq \dd \}$ for some $\dd$ \cite[Theorem
17.4]{Sch}. The set of rows of $M$ can be chosen to be
$$\{ \mm \in \ZZ^n \,:\, \mm = \yy A, \, \yy \geq {\bf 0}, \, ||\mm||_{\infty}
\leq n^{2n} \Delta^n \}$$ where $\Delta$ is the maximum absolute value
of a minor of $A$. In practice, $M$ could be much smaller.  For
instance if $A$ is the $4 \times 2$ matrix with rows $(1, 2)$,
$(-2,-3)$, $(1,0)$ and $(0,1)$, it suffices to augment $A$ with the
rows $(1,1),(0,-1),(-1,-2),(-1, -1)$, while $n^{2n} \Delta^n = 144$.

In Section~\ref{sec:definitions} we introduce a vector version of the
Chv{\'a}tal-Gomory procedure called {\em iterated basis normalization}
(IBN) that constructs a sufficient $M$ from the matrix $A$. The {\em
  small Chv{\'a}tal rank} (SCR) of $A$ is the number of rounds of IBN
necessary to generate this $M$. A similar definition can be made when
$\bb$ is fixed. The SCR of $A$ (respectively of $A{\bf x} \leq \bb$)
is at most its Chv{\'a}tal rank even though IBN may not terminate when
$n \geq 3$. We show that in every dimension, there are systems $A \xx
\leq \bb$ for which SCR is two while the Chv{\'a}tal rank is
arbitrarily high.

In Section~\ref{sec:supernormal} we completely characterize matrices
for which SCR is zero.  This requires the notion of {\em
  supernormality} introduced in \cite{HMS} which generalizes the
familiar notion of unimodularity.  We produce a family of matrices of
increasing dimension for which SCR is zero but Chv{\'a}tal rank is not
zero.

In Section~\ref{sec:stable set problem} we apply the theory of SCR to
$\Frac(G)$, the fractional stable set polytope of a graph $G$. We
determine the structure of the vectors produced by IBN in rounds one
and two. As a consequence we see that the normals of many of the
well-known facet inequalities of the stable set polytope, $\stab(G)$,
appear within two rounds of IBN. It is a long-standing open problem to
describe $\stab(G)$ when $G$ is a claw-free graph. We show that many
of the complicated facet normals of $\stab(G)$ when $G$ is claw-free
appear in two rounds of IBN which reveals a uniform hypercyclic
structure in these ad hoc examples.

Section~\ref{sec:lowerbounds} contains lower bounds for SCR which
contrast with the results in the earlier sections.  We show that if $n
\geq 3$, SCR may grow exponentially in the bit size of the matrix $A$,
asymptotically just as fast as Chv{\'a}tal rank. For polytopes in the
unit cube $[0,1]^n$, SCR can be at least $n/2$. We also exhibit a
lower bound that depends on $n$ for $\scr(\Frac(G))$ as $G$ varies
over all graphs with $n$ vertices. A brief discussion of possible
upper bounds and computational evidence supporting our guesses are
also provided.


The SCR of $A \xx \leq \bb$ or $A$ offers a coarser measure than
Chv{\'a}tal rank of the complexity of the integer programs associated
to them. The goal here is to determine how quickly the facet normals
of an integer hull are produced from the normals of the rational
polyhedron, ignoring the right-hand-sides of the facet
inequalities. Our main message is that, in many cases, facet normals
are produced surprisingly fast by the Chv{\'a}tal-Gomory procedure but
the right-hand-side can take a long time to be computed, which makes
Chv{\'a}tal rank high. The coarseness of SCR can be a powerful
organizational tool that can reveal the unifying structure behind
seemingly ad hoc facet normals of a class of examples. An illustration
of this philosophy can be found in Example~\ref{ex:clawfree} where we
show that many difficult facet normals that have been found for the
stable set polytope of a claw-free graph are produced within two
rounds of IBN. While the Chv{\'a}tal-Gomory procedure carries along
both the number theoretic and geometric parts of an integer hull
computation, SCR focuses on the number theory alone, often
revealing interesting structural facts that are difficult to see through
the fine Chv{\'a}tal-Gomory lens.

\section{Main Definitions.} \label{sec:definitions} 

Fix a matrix $A \in \ZZ^{m \times n}$ of rank $n$ and let $\A =
\{\aaa_1, \ldots, \aaa_m\}$ \ be the vector configuration in $\ZZ^n$
consisting of the rows of $A$. We assume that each row of $A$ is
primitive (i.e., the gcd of its components is one). For each $\bb \in
\ZZ^m$, consider the rational polyhedron $Q_\bb := \{ \xx \in \RR^n
\,:\, A \xx \leq \bb \}$ and its integer hull $Q_\bb^I :=
\textup{conv}(Q_\bb \cap \ZZ^n)$ where $\textup{conv}$ denotes convex
hull.  Since $\textup{rank}(A) = n$, every minimal face of $Q_\bb$,
and $Q_\bb^I$ (if non-empty), is a vertex.  A {\em Hilbert basis} of a
rational polyhedral cone $K \subseteq \RR^n$ is a set $\{\hh_1,\ldots,
\hh_t\} \subset K \cap \ZZ^n$ such that if $\kk \in K \cap \ZZ^n$ then
$\kk = \sum_{i=1}^{t} n_i \hh_i$ where $n_i \in \NN$. If $K$ is
pointed then it has a unique minimal Hilbert basis.  Write $\Hilb(K)$
(respectively, $\Hilb(\A)$) for a minimal Hilbert basis of $K$ (respectively, 
$\cone(\A)$).

The Chv{\'a}tal-Gomory procedure \cite{Chv73}, \cite[\S 23]{Sch} for
computing $Q_\bb^I$ works as follows. For each vertex $\vv$ of
$Q_{\bb}$, set $\A_{\vv} := \{ \aaa_i \in \A: \aaa_i \cdot \vv = b_i \}$
and define $Q_{\bb}^{(1)}$ to be the polyhedron cut out by the
inequalities $\hh \cdot \xx \leq \lfloor \hh \cdot \vv \rfloor$ for
every vertex $\vv$ of $Q_{\bb}$ and every vector $\hh \in
\Hilb(\A_\vv)$. Then $Q_{\bb}^I \subseteq Q_{\bb}^{(1)} \subseteq
Q_{\bb}$.  For $i \geq 2$, define $Q_{\bb}^{(i)} :=
(Q_{\bb}^{(i-1)})^{(1)}$. For a positive integer $k$, $Q_\bb^{(k)}$ is
called the $k$-th {\em Chv{\'a}tal closure} of $Q_\bb$.  The {\em
  Chv{\'a}tal rank} of $A \xx \leq \bb$ (equivalently, $Q_\bb$) is the
smallest number $t$ such that $Q_{\bb}^{(t)} = Q_\bb^I$. This rank
only depends on $Q_\bb$ and not the inequality system defining it. The
{\em Chv{\'a}tal rank} of $A$ is the maximum over all $\bb \in \ZZ^m$
of the Chv{\'a}tal ranks of $A \xx \leq \bb$. The Chv{\'a}tal-Gomory
procedure and the Chv{\'a}tal ranks are all finite \cite[Chapter
23]{Sch}.

To study just the facet normals of the integer hulls
$Q_\bb^I$ for every $\bb$, we modify the Chv{\'a}tal-Gomory procedure as
follows. An $n$-subset $\tau \subseteq [m] := \{1,2,\ldots,m\}$ is
called a {\em basis} if the submatrix $A_\tau$, consisting of the rows
of $A$ indexed by $\tau$, is non-singular. Let $\A_\tau$ be the set of
rows of $A_\tau$.  We call $\cone(\A_\tau)$ a {\em basis cone} since
$\A_\tau$ is a basis of $\RR^n$. The set $\A$ contains at least one
basis cone since $\textup{rank}(A) = n$.

\begin{observation} \label{obs:unionofhbs}
  Suppose $\sigma \subseteq [m]$ such that $\A_\sigma$
  linearly spans $\RR^n$. Then the union of the minimal Hilbert bases
  of the basis cones $\cone(\A_\tau)$, as $\tau$ varies over the bases
  contained in $\sigma$, is a Hilbert basis for $\cone(\A_\sigma)$.
\end{observation}

\begin{algorithm} {\em Iterated Basis Normalization
    (IBN)}\label{alg:IBN}\\
\noindent{\em Input}: $A \in \ZZ^{m \times n}$ satisfying the
assumptions above. 

\begin{enumerate} 
\item Set $\A^{(0)} := \A$.
\item For $k \geq 1$, let $\A^{(k)}$ be the union of all the (unique)
minimal Hilbert bases of all basis cones in $\A^{(k-1)}$.
\item If $\A^{(k)} = \A^{(k-1)}$, then stop. Otherwise repeat. 
\end{enumerate}
\end{algorithm}

\begin{remark}
Since each vector in $\A$ is primitive, $\A \subseteq \A^{(1)}.$ 
Every vector created during IBN is also primitive and so 
$\A \subseteq \A^{(1)} \subseteq \A^{(2)} \subseteq \ldots$.
\end{remark}

\begin{lemma} \label{lem:all non-negative} If all elements of $\A$ are
  non-negative except for the negative unit vectors $-\ee_i, i \in
  [n]$, then for each non-negative integer $k$, all vectors in
  $\A^{(k)}$ besides the original $-\ee_i$'s are also non-negative.
\end{lemma}

\begin{proof}
  The claim holds for $k=0$, and suppose it holds up to $k-1$.
  When IBN constructs $\A^{(k)}$ from $\A^{(k-1)}$, for each
  $i \in [n]$, the only vector available with negative $i$-th
  coordinate is $-\ee_i$ but since its multiplier lies in $[0,1)$, the
  $i$th coordinate of the resulting Hilbert basis elements cannot be
  negative.
 \end{proof}

Let $A^{(k)}$ denote a matrix whose rows are the elements of
$\A^{(k)}$ with the rows in $\A^{(k)} \setminus \A^{(k-1)}$ appended
at the bottom of $A^{(k-1)}$.

\begin{definition} \label{def:SCR} $ $
\begin{enumerate}
\item The {\em small Chv{\'a}tal rank (SCR)} of the system of 
inequalities $A \xx \leq \bb$ defining $Q_\bb$ is the smallest number $k$ 
such that there is an integer vector $\bb^\prime$ satisfying 
$$Q_\bb^I \, = \, \{\xx \in \RR^n \, : \, A^{(k)} \xx \leq \bb^\prime \}.$$
\item The SCR of a matrix $A$ is the supremum of the SCRs of all
  systems of the form $A \xx \leq \bb$ as $\bb$ varies in $\ZZ^m$.
\end{enumerate}
\end{definition}

\begin{proposition} \label{prop:ranks}
For any $\bb \in \ZZ^m$, the SCR of $A \xx \leq \bb$ is at most the
Chv{\'a}tal rank of the same system, and the SCR of $A \in \ZZ^{m
\times n}$ is at most the Chv{\'a}tal rank of $A$.  In particular,
the SCR is always finite.
\end{proposition}

\begin{proof}
  If $\vv$ is a vertex of some intermediate polyhedron $Q_\bb^{(i)} =
  \{ \xx \,:\, U \xx \leq \uu \}$ in the Chv{\'a}tal-Gomory procedure,
  then $\U_\vv$ linearly spans $\RR^n$. By Observation~\ref{obs:unionofhbs}
  and induction, a Hilbert basis of $\U_\vv$ is contained in
  $\A^{(i+1)}$ and therefore, $\A^{(i+1)}$ contains the normals of an
  inequality system describing $Q_\bb^{(i)}$. In particular, if the
  Chv{\'a}tal rank of $A \xx \leq \bb$ is $t$, then the normals of an
  inequality system describing $Q_{\bb}^I$ are in $\A^{(t)}$. 
 \end{proof}

\begin{lemma} \label{lem:ibn-codim2}
When $n=2$, $\A^{(2)} = \A^{(1)}$, and IBN terminates in one round.
\end{lemma}

\begin{proof} Pick $\rr, \sss \in \A^{(1)} \subset \ZZ^2$ such that 
  $\cone(\rr,\sss)$ is a basis cone. Let 
$$ \ttt_1:=\rr,\ttt_2,\ldots,\ttt_{k-1},\ttt_k:=\sss $$ 
be the elements of $\A^{(1)}$ in
  $\cone(\rr,\sss)$ in cyclic order from $\rr$ to $\sss$. Then for each $i \in
  \{1,\ldots,k-1\}$, $\cone(\ttt_i,\ttt_{i+1})$ is unimodular. (This is an
  artifact of $\RR^2$. See
  \cite[Corollary 3.11]{MT02} for a proof.)  Hence a Hilbert basis of
  $\cone(\rr,\sss)$ is contained in $\{ \ttt_1, \ldots, \ttt_k \}$, and
  $A^{(2)} = \A^{(1)}$.
 \end{proof}

\begin{corollary}\label{cor:scr-codim2}
If $A \in \ZZ^{m \times 2}$, then the SCR of $A$ is at most one.
\end{corollary}

\begin{example}
In contrast, Chv{\'a}tal rank can be arbitrarily large even for 
$A \in \ZZ^{3 \times 2}$. Fix $j \geq 1$ and consider the system $A
\xx \leq \bb$ where 
$$A = \left( \begin{array}{rr}  -1 & 0 \\ 1 & 2j \\ 1 & -2j
    \end{array} \right) \, \textup{and} \; \bb = (0, 2j, 0)^t. $$  
The polyhedron $Q_\bb$ is a triangle in $\RR^2$ with vertices $(0,0)$,
$(0,1)$ and $(j,1/2)$, and $Q_\bb^I$ is the line segment from $(0,0)$
to $(0,1)$. It is noted in \cite[\S 23.3]{Sch} that 
the Chv{\'a}tal rank of $A \xx \leq \bb$ is at least $j$. 

Fix $n \geq 2$ and $m \geq n+1$.  By taking the product of $Q_\bb$
from above with the $(n-2)$-dimensional positive orthant and then
adjoining $m-n-1$ redundant inequalities, we can produce $A' \xx \leq
\bb'$, $A' \in \ZZ^{m \times n}$ with the same property that SCR is
one but Chv\'atal rank is arbitrarily large.
\end{example}

Unlike for $n \leq 2$, IBN need not terminate when $n > 2$.

\begin{example}
  Take $\A = \{(0,3,1),(1,1,1),(2,5,5),(1,4,3)\}$.  For each positive  
  integer $k$, set
$$\uu_k := (k, 2k+2, 2k+1)\,\,\, \textup{and} \,\,\, \vv_k := (k,
2k+1, 2k).$$
Note that $\uu_1 = (1,4,3)$ is a row of $A$.  To show
that IBN does not terminate on $\A$, one can check the following two
assertions. We omit the details.

\begin{enumerate}
\item For each $k \geq 1$, $\vv_k \in \Hilb(\{(0,3,1), (1,1,1),
  \uu_k\})$.
\item For each $k \geq 1$, $\uu_{k+1} \in \Hilb(\{(0,3,1), (2,5,5),
  \vv_k \})$.
\end{enumerate}
A second such example appears in \cite{HMS}.
\end{example}

Despite this example, the SCR of any matrix or system of inequalities is finite, and we will illustrate ways to bound it in many instances. 

\begin{definition} \label{def:small chvatal closure} For a positive
  integer $k$, the $k$-th {\bf small Chv{\'a}tal closure} of $A
  \xx \leq \bb$ is the set $S_{A,\bb}^{(k)} := \{ \xx \in \RR^n \,:\, \aaa
  \cdot \xx \leq m_\aaa, \,\,\aaa \in \A^{(k)} \}$ where $m_\aaa :=
  \textup{max} \{ \aaa \cdot \xx \,:\, A \xx \leq \bb, \,\,\xx \in
  \ZZ^n \}$. 
\end{definition}

This is a definition for inequality systems: if $P = \{\xx \in \RR^n
\,:\, A \xx \leq \bb\} = \{ \xx \in \RR^n \,:\, A' \xx \leq \bb' \}$,
then for a given $k$, $S_{A,\bb}^{(k)}$ may not equal
$S_{A',\bb'}^{(k)}$. However, for a fixed $A \xx \leq \bb$,
$S_{A,\bb}^{(k)} \subseteq Q_\bb^{(k)}$ for each non-negative integer
$k$.

\begin{lemma} \label{lem:poly characterization of scr} The SCR of $A
  \xx \leq \bb$ is the smallest integer $k$ such that $Q_\bb^I =
  S_{A,\bb}^{(k)}$.
\end{lemma}

\begin{proof} If $\scr(A \xx \leq \bb) = k$, then $Q_\bb^I = \{ \xx
  \in \RR^n \,:\, \aaa \cdot \xx \leq b_\aaa, \,\,\aaa \in \A^{(k)} \}$
  for some scalars $b_\aaa$. However, $b_\aaa \geq m_\aaa$ for each $\aaa
  \in \A^{(k)}$, so $S_{A,\bb}^{(k)} \subseteq \{ \xx \in
  \RR^n \,:\, A \xx \leq \bb \}^I$, and hence they are equal. On the
  other hand, $S_{A,\bb}^{(k-1)} \neq Q_\bb^I$ since otherwise $\scr(A
  \xx \leq \bb)$ would be less than $k$. 
 \end{proof}

\section{Matrices with small Chv{\'a}tal rank
  zero}\label{sec:supernormal}  

We begin our study of SCR by characterizing the matrices $A$ for which
SCR is zero. These are precisely the $A$'s with the property that for
each $\bb \in \ZZ^m$, there is a $\bb' \in \ZZ^m$ such that $Q_\bb^I =
\{ \xx \in \RR^n \,:\, A \xx \leq \bb' \}$. Our characterization
offers a generalization of the familiar notion of unimodularity.

\begin{definition}\label{def:unimodular}
  A vector configuration $\A$ in $\ZZ^n$ is {\em unimodular} if for
  every subset $\A'$ of $\A$, $\A'$ is a Hilbert basis for
  $\textup{cone}(\A')$.
\end{definition}


\begin{definition} \label{def:tdi} \cite[Theorem~22.5]{Sch}
  A system of linear inequalities $A \xx \leq \bb$ is {\em totally
    dual integral (TDI)} if the set
$ \A_F := \{ \aaa_i \in \A \,:\, \aaa_i \xx  = b_i \,\, \forall \,\, \xx \in F \}$ is a Hilbert basis of the cone it generates for every face $F$ of the polyhedron $Q_{\bb} = \{ \xx \,:\, A \xx \leq \bb \}$.
\end{definition}

The following characterizations of matrices with Chv{\'a}tal rank zero
are well-known, while characterizations of higher Chv{\'a}tal rank are
unknown.

\vspace{.2cm}

\begin{theorem} \label{thm:unimodular} \cite{Sch} 
  Let $\A = \{ \aaa_1, \ldots, \aaa_m \} \subset \ZZ^n$ be such that the
  matrix $A$ whose rows are $\aaa_1, \ldots, \aaa_m$ has rank $n$. Then
  the following are equivalent:
\begin{enumerate}
\item $\A$ is unimodular.
\item Every basis in $\A$ is a basis of $\ZZ^n$ as a lattice.
\item Every (regular) triangulation of $\A$ is unimodular.
\item For all $\bb \in \ZZ^m$, the inequality system $A \xx \leq \bb$
  is TDI. 
\item For all $\bb \in \ZZ^m$, the polyhedron $Q_{\bb} = \{ \xx \in
  \RR^n \,:\, A \xx \leq \bb \}$ is integral.
\item The Chv{\'a}tal rank of $A$ is zero.
\end{enumerate}
\end{theorem}

Theorem~\ref{thm:supernormal} will provide a complete analogue to
Theorem~\ref{thm:unimodular} when SCR replaces Chv{\'a}tal rank. A
vector configuration $\A$ in $\ZZ^n$ is {\em normal} if it is a
Hilbert basis for $\cone(\A)$. 

\begin{definition} \cite{HMS} \label{def:supernormal}
  A configuration $\A$ is {\em supernormal} if for every subset
  $\A^\prime$ of $\A$, $\A \cap \cone(\A^\prime)$ is a Hilbert basis
  of $\cone(\A^\prime)$.
\end{definition}

Following \cite{HMS}, we say that a system $A \xx \leq \bb$ is {\em
  tight} if for each $i = 1, \ldots, m$, the hyperplane $\aaa_i \xx =
b_i$ contains an integer point in $Q_\bb$ and hence supports
$Q_{\bb}^I$. When the inequality system is clear, we
simply say that the polyhedron $Q_{\bb} = \{ \xx \,:\, A \xx \leq \bb
\}$ is tight. If $Q_\bb^I$ is nonempty, recall that 
$$m_{\aaa_i} := \textup{max} \{ \aaa_i \cdot \xx \,:\, \xx \in Q_{\bb}^I
\} \,\,\,\textup{for} \,\,\, i = 1, \ldots, m$$ and set $\beta :=
(m_{\aaa_i}) \in \ZZ^m$. Then $Q_\bb \supseteq Q_\beta \supseteq Q_\bb^I$
and $Q_\beta$ is tight.
  
\begin{theorem} \label{thm:supernormal}
  Let $\A = \{ \aaa_1, \ldots, \aaa_m \} \subset \ZZ^n$ be
  a configuration of primitive vectors such that the matrix $A$ whose 
  rows are $\aaa_1, \ldots, \aaa_m$ has rank $n$. Then the following are 
  equivalent.
\begin{enumerate}
\item $\A$ is supernormal.
\item Every basis $\A'$ in $\A$ has the property that  
  $\A \cap \textup{cone}(\A')$ is a Hilbert basis of $\textup{cone}(\A')$,
  or equivalently, $\A = \A^{(1)}$.
\item Every (regular) triangulation of $\A$ that uses all the vectors
  is unimodular.
\item For all $\bb \in \ZZ^m$, $A \xx \leq \bb$ is TDI whenever
  $Q_{\bb}$ is tight.
\item For all $\bb \in \ZZ^m$, the polyhedron $Q_{\bb}$ is integral whenever
$Q_\bb$ is tight. 
\item The SCR of $A$ is zero.
\end{enumerate}
\end{theorem}

The equivalence of (1), (3), and (4) is shown in~\cite[Proposition 3.1
and Theorem 3.6]{HMS}. Our contribution is the remaining set of
equivalences.

\begin{proof}

[(1) $\Rightarrow$ (2)]: This is immediate from the definition of 
supernormality. 

[(2) $\Rightarrow$ (3)]: Let $T$ be a triangulation of $\A$ using all
of the vectors and $\sigma$ index a maximal simplex of $T$.  Then the
sub-configuration $\A_\sigma$ is a basis of $\A$ and by (2), $\A$
contains $\Hilb(\A_\sigma)$.  But since every vector in $\A$ is used
in the triangulation $T$, none can lie inside or on the boundary of
$\cone(\A_\sigma)$ except those in $\A_\sigma$ itself.  Thus
$\A_\sigma$ is the Hilbert basis of its own cone.  This implies that
$\A_\sigma$ is a lattice basis, so $\sigma$ is a unimodular simplex.
Since $\sigma$ was arbitrary, $T$ is a unimodular triangulation.

[(4) $\Rightarrow$ (5)]: This follows from~\cite[Corollary 22.1c]{Sch},
which says that for a  $\bb \in \ZZ^m$, if $A \xx \leq \bb$ is TDI,
then $Q_\bb$ is integral. 

[(5) $\Leftrightarrow$ (6)]: Suppose $Q_\bb$ with $\bb \in \ZZ^m$ is
integral whenever it is tight. Then for $\bb \in \ZZ^m$ with
$Q_{\bb}^I \neq \emptyset$, $Q_{\beta}$ is integral since it is tight.
But $Q_{\bb}^I \subseteq Q_{\beta} \subseteq Q_{\bb}$ which implies
that $Q_{\beta} = Q_{\bb}^I$ and the SCR of $A$ is zero.

Suppose the SCR of $A$ is zero and some $Q_{\bb}$ is tight.  Then no
new facet normals are needed for $Q_{\bb}^I$, so $Q_\bb^I = Q_{\beta}
\subseteq Q_{\bb}$. Since for $i=1, \ldots, m$, $\aaa_i \cdot \xx =
b_i$ and $\aaa_i \cdot \xx = m_{\aaa_i}$ both support $Q_{\bb}^I$,
$\beta = \bb$. Thus $Q_\bb = Q_\beta$ is integral.

[(6) $\Rightarrow$ (3)]: Suppose there exists a non-unimodular
(regular) triangulation $T$ of $\A$ that uses all the vectors in $\A$.
Let $\A'$ be a basis in $\A$ whose elements form a non-unimodular
facet in $T$ and let $A'$ be the non-singular square matrix whose rows
are the elements of $\A'$. Then there exists a $\bb' \in \ZZ^n$ such
that $\{\xx \in \RR^n \,:\, A' \xx \leq \bb'\}$ is tight and its
unique vertex is not integral. Since no element of $\A \backslash \A'$
lies in $\textup{cone}(\A')$, by choosing very large right-hand-sides
for the elements in $\A \backslash \A'$, one gets a $Q_\bb$ in which
the fractional vertex of $\{ \xx \in \RR^n \,:\, A' \xx \leq \bb'\}$
and its neighborhood survive. Further, $\bb$ can be chosen so that
$Q_\bb$ is tight. Therefore, the SCR of $A$ is not zero.
 \end{proof}

\begin{example} If the rows of $A$ are not primitive 
  then supernormality is not necessary for the SCR of $A$ to be zero.
  Take $\A = \{(2,0),(0,2)\}$.  Then for each $\bb \in \ZZ^2$,
  $Q_{\bb} = \{ (x_1,x_2) \in \RR^2 \,:\, x_1 \leq \frac{b_1}{2}, \,\,
  x_2 \leq \frac{b_2}{2} \}$. Hence $Q_{\bb}$ is tight if and only if
  both $b_1$ and $b_2$ are even, in which case it has the unique
  integer vertex $(\frac{b_1}{2}, \frac{b_2}{2})$. Therefore all tight
  $Q_{\bb}$'s are integral but $A$ is not supernormal. It is easy to
  see that the SCR of $A$ is zero.
\end{example}

\begin{remark} \label{rem:complexity}
  If the dimension $n$ is fixed, then it is possible to determine
  whether $\A$ is supernormal (and hence whether SCR is zero) in
  polynomial time. The number of basis cones is at most
  $m \choose n$, so it suffices to check whether $\A \cap
  \textup{cone}(\A')$ is normal for each basis $\A'$ in $\A$. Barvinok
  and Woods~\cite[Theorem 7.1]{BaWo} show that in fixed dimension, a
  rational generating function for the Hilbert basis of each cone can
  be computed in polynomial time. We then subtract the polynomial
  $\sum_{\aaa \in \cone(\A') \cap \A} \xx^\aaa$ from this rational
  function, square the difference, and evaluate at $\xx =
  (1,\ldots,1)$.  This can also be done in polynomial time
  ~\cite[Theorem 2.6]{BaWo} and the result is zero if and only if $\A
  \cap \textup{cone}(\A')$ is normal.
\end{remark}

\begin{problem} Can one recognize the supernormality of $A$ in
  polynomial time analogous to Seymour's result for unimodularity
  \cite{Seymour}, \cite[Chapter 20]{Sch}?  
\end{problem} 

We close this section with a family of matrices for which Chv{\'a}tal
rank is not zero while SCR is. The existence of such families was a
question in \cite{HMS}.

\begin{proposition} \label{prop:supernormalnotunimodular}  
There exist configurations in arbitrary dimension which are
supernormal but not unimodular.  
\end{proposition}

\begin{proof}
Let $k$ be a positive integer and $\A$ be the rows of the 
$(2k+1) \times (2k+1)$ matrix
$$
A = \left( \begin{array}{rrrrrr}
    1 & 1 & 0 & \ldots & 0 & 0 \\
    0 & 1 & 1 & \ldots & 0 & 0 \\
    &   &   & \vdots &   &   \\
    0 & 0 & 0 & \ldots & 1 & 1 \\
    1 & 0 & 0 & \ldots & 0 & 1 \end{array} \right).$$
That is, $A$ is
the edge-vertex incidence matrix of an odd circuit.  The determinant of
$A$ is two, so there is exactly one Hilbert basis element of
$\cone(\A)$ that does not generate an extreme ray: the all-ones
vector {\bf 1}.

We claim that all maximal minors of $A^{(1)}$ except for
$\textup{det}(A)$ are $\pm 1$.  This implies that $\A^{(1)}$ equals
$\A^{(2)}$, and hence by Theorem~\ref{thm:supernormal}, $\A^{(1)}$
is supernormal.  But since $\A$ is not unimodular, neither is
$\A^{(1)}$, proving the proposition.

To prove the claim, by symmetry it suffices to check a single minor of
$A^{(1)}$ different from $\textup{det}(A)$, for instance the minor
obtained by removing the last row of $A$ from $A^{(1)}$.  By cofactor
expansion on the last column, this minor equals 
$\textup{det}(D_2) - \textup{det}(D_1)$ 
where $D_1$ and $D_2$ are the $2k \times 2k$ matrices
$$ 
\begin{array}{clc}
D_1 = \left( \begin{array}{rrrrrr} 
1 & 1 & 0 & \ldots & 0 & 0 \\
0 & 1 & 1 & \ldots & 0 & 0 \\
  &   &   & \vdots &   &   \\
0 & 0 & 0 & \ldots & 1 & 1 \\
1 & 1 & 1 & \ldots & 1 & 1 \\ \end{array} \right) &
\textup{and} &
D_2 = \left( \begin{array}{rrrrrr}
    1 & 1 & 0 & \ldots & 0 & 0 \\
    0 & 1 & 1 & \ldots & 0 & 0 \\
    &   &   & \vdots &   &   \\
    0 & 0 & 0 & \ldots & 1 & 1 \\
    0 & 0 & 0 & \ldots & 0 & 1 \\ \end{array} \right).
\end{array} $$
The last row of $D_1$ is the sum of its odd-indexed rows
so $\textup{det}(D_1) = 0$. Further, $D_2$ is upper triangular
with 1's on the diagonal, so $\textup{det}(D_2) = 1$.
 \end{proof}


\section{Application to the stable set problem in a
  graph} \label{sec:stable set problem} 

We now apply the theory of SCR in the specific context of the {\em
  maximum stable set problem} in a graph. Besides being an important
example, the results offer a glimpse of the kind of insights that
might be possible when SCR is examined for problems with structure. We
will show that the normals of many well-known valid inequalities of
the stable set polytope appear within two rounds of IBN.

Let $G=([n],E)$ be an undirected graph with vertex set $[n]$ and edge
set $E$. A {\em stable set} in $G$ is a subset $U \subseteq [n]$ such
that $\{i,j\} \not \in E$ for any pair $i,j \in U$. The {\em stability
  number} $\alpha(G)$ is the maximum size of a stable set in $G$, and
the {\em stable set problem} seeks a stable set in $G$ of cardinality
$\alpha(G)$. This is a well-studied, NP-hard problem in combinatorial
optimization that has been approached via linear and semidefinite
programming. The basic idea behind both approaches is as follows. Let
$\ee_i$ denote the $i$th standard unit vector in $\RR^n$ and $\ee(U)
:= \sum_{i \in U} \ee_i \in \{0,1\}^n$ be the {\em characteristic
  vector} of $U \subseteq [n]$. The convex hull of the characteristic
vectors of all stable sets in $G$ is the {\em stable set polytope},
$\stab(G)$, and the stable set problem can be modeled as the linear
program:
\begin{equation}\label{eqn:maxstablesetprob}
  \alpha(G) = \textup{max}\{ \sum_{i=1}^{n} x_i  \,:\, \xx \in \stab(G) \}.
\end{equation}
The polytope $\stab(G)$ is not known a priori, and so the linear and
semidefinite programming approaches construct successive outer
approximations of $\stab(G)$ that eventually yield an optimal solution
of (\ref{eqn:maxstablesetprob}). The linear programming relaxations of
$\stab(G)$ are all polytopes and the standard starting approximation
is the {\em fractional stable set polytope}
$$\Frac(G) := \{ \xx \in \RR^n \,:\, x_i + x_j
\leq 1 \,\,\,(\forall \{i,j\} \in E), \,\, x_i \geq 0 \,\,\,(\forall i
\in [n]) \}$$ whose integer hull is $\stab(G)$. See \cite[Chapter
9]{GLS} for more details.

In this section we examine the SCR of the inequality system defining
$\Frac(G)$ which we denote as $\scr(\Frac(G))$ since the inequality
system is well defined. The input to IBN is $$\A_G := \{ \ee_i + \ee_j
\,:\, \{i,j\} \in E \} \cup \{ - \ee_i \,:\, i \in [n] \},$$ and let
$\A_G^{(k)}$ be the configuration created by IBN after $k$ rounds. We
will describe $\A_G^{(1)}$ and $\A_G^{(2)}$ combinatorially and show
that $\A_G^{(2)}$ contains the normals of many well-known classes of
facet inequalities of $\stab(G)$.

For $U \subseteq [n]$, let $\xx(U) := \sum_{i \in U} x_i$. If $H =
(V_H,E_H)$ is a subgraph in $G$ then we write $\xx(H)$ for $\xx(V_H)$
and $\ee(H)$ for $\ee(V_H)$. By a {\em circuit} in $G$ we mean a cycle
(closed walk) in $G$ with distinct vertices and edges. A {\em hole} in
$G$ is a chordless circuit and an {\em antihole} is the complement of
a hole. A {\em wheel} in $G$ is a circuit with an additional vertex
$u_0$ that is joined by edges to all vertices of the cycle. The wheel
$W$ is odd if $|V_W \setminus \{u_0\}|$ is odd. The following are
well-known classes of valid inequalities of $\stab(G)$:
\begin{center}
\begin{tabular}{lll}
  1. {\bf non-negativity}   &  $x_i \geq 0$, $i \in [n]$ \\
  2. {\bf edge}  & $x_i + x_j \leq 1$,\,\,\, $\{i,j\} \in E$ \\
  3. {\bf clique}  & $\xx(K) \leq 1$,\,\,\, $K$ clique in
  $G$\\ 
  4. {\bf odd hole/circuit}  & $\xx(C) \leq \lfloor \frac{|C|}{2}
  \rfloor$,\,\,\, $C$ odd hole/circuit in $G$\\
  5. {\bf odd antihole}  & $\xx(A) \leq 2$,\,\,\, $A$ an odd antihole
  in $G$\\
  6. {\bf rank}  & $\xx(H) \leq \alpha(H)$,\,\,\, $H$ a subgraph in $G$\\
  7. {\bf odd wheel}  & $\xx(V_W \setminus \{u_0\}) +
  \frac{|V_W|-2}{2} x_{u_0} \leq \frac{|V_W|-2}{2}$, $W$ a wheel in $G$.\\
\end{tabular}
\end{center}

\vspace{.2cm}

Constraints 1-5 are all rank inequalities while the odd wheel
inequalities are not.  Our interest will be in determining the least
$k$ for which the normal of a valid inequality for $\stab(G)$ appears
in $\A_G^{(k)}$. 

For a graph $G$, let $Q_G^{(k)}$ denote the $k$-th Chv{\'a}tal closure
of $\Frac(G)$ and $S_G^{(k)}$ denote the $k$-th small Chv{\'a}tal
closure of the inequality system defining $\Frac(G)$. Then
$$\Frac(G) = Q_G^{(0)} = S_G^{(0)} \,\,\,\textup{and} \,\,\, Q_G^{(k)}
\supseteq S_G^{(k)} \,\,\,\forall \,\,\,k.$$ Note that $S_G^{(0)}$ is
obtained by making the inequality system defining $\Frac(G)$ {\em
  tight} in the sense of Section~\ref{sec:supernormal}, but this
inequality system is already tight, and so $Q_G^{(0)} = S_G^{(0)}$.  We
now determine the structure of $\A_G^{(1)}$.

\begin{proposition} \label{prop:A1} The elements of $\A_G^{(1)}
  \setminus \A_G^{(0)}$ are precisely the characteristic vectors,
  $\ee(C)$, of odd circuits $C$ in $G$.
\end{proposition}

\begin{proof} Suppose $\vv \in \A_G^{(1)} \setminus \A_G^{(0)}$. Then
  there is a basis $\B = \{\bb_1, \ldots, \bb_n\} \subseteq
  \A_G^{(0)}$ such that $\vv = \sum_{i=1}^{n}\lambda_i \bb_i$ with $0
  \leq \lambda_i < 1$. Let 
  $$\sigma := \{ i \in [n]: \lambda_i > 0 \}, \, W := \{ j \in [n] \,:\, (\bb_i)_j \neq 0 \textup{ for some } i \in \sigma \},$$ 
and $B'$ be the submatrix of $B$ whose rows are indexed
  by $\sigma$ and columns by $W$. (Recall that $B$ is the matrix with
  rows $\bb_1, \ldots, \bb_n$.) Then $ \textup{span} \{\bb_i \,: i \in
  \sigma \} \subseteq \textup{span} \{\ee_j \, : j \in W \} $ and
  since $\{\bb_i \,:\, i \in \sigma \}$ is part of a basis, we obtain
  $ \left| \sigma \right| \leq \left| W \right|$. Also, since $\vv \in
  \ZZ^{n}$ and $\lambda_i \notin \ZZ$ for every $i \in \sigma$, for
  every $j \in W$ there must be at least two rows in $B'$ whose $j$th
  entries are nonzero. However, each $\aaa_i \in \A_G$ has at most two
  nonzero coordinates. Thus if $k$ is the total number of nonzero
  entries in $B'$, we have
$$ 2 \left| W \right| \leq k \leq 2 \left| \sigma \right| \leq 2
\left| W \right| $$  
and so each inequality must be satisfied with equality. This means
that:

\begin{enumerate}
\item for every $i \in \sigma$, $\aaa_i$ has exactly two nonzero
  entries: it is the incidence vector of an edge in $G$; and
\item for every $j \in W$, the column of $B'$ indexed by $j$ has
  exactly two nonzero entries.
\end{enumerate}

That is, $B'$ is the incidence matrix of a subgraph of $G$ in which
every vertex has degree two, and so it is a union of disjoint circuits
in $G$.  If any of these circuits is even then the corresponding rows
of $B_\sigma$ are dependent, which is a contradiction. Also, if there
is more than one odd circuit, then $\vv$ is the sum of at least two
different integer vectors in the fundamental parallelepiped spanned by
the rows of $B_\sigma$, which contradicts that $\vv$ is in
$\Hilb(\B_\sigma)$. Therefore, there is a single odd circuit in $G$
with vertex set $W$. It is now a simple exercise to see that $\vv$ is
the all ones vector and that for all $i \in \sigma$, $\lambda_i =
1/2$.

Conversely, if $C$ is an odd circuit in $G$, then by taking $E_\B$ to
be the collection of edges in $C$ and augmenting the corresponding
elements of $\A_G$ to a basis by adding $-\ee_i$'s indexed by vertices
outside $C$, we produce $\ee(C)$ as an element of $\A_G^{(1)}$.
 \end{proof}

\begin{corollary} \label{cor:all non-negative} For each non-negative
  integer $k$, all vectors in $\A_G^{(k)}$ different from the
  $-\ee_i$'s are non-negative.
\end{corollary}

\begin{proof} Follows from Lemma~\ref{lem:all non-negative}.
 \end{proof}

\begin{corollary} \label{cor:oddcircuit normal} Let $H$ be an induced
  subgraph of $G$ such that there is an odd circuit in $G$ through the
  vertices of $H$.  Then $\ee(H)$, the normal of the rank inequality $
  \xx(H) \leq \alpha(H)$, appears in $\A_G^{(1)}$. In particular, the
  normals of all odd hole and odd clique inequalities appear in
  $\A_G^{(1)}$.
\end{corollary}

\begin{corollary} \label{cor:first small chvatal closure} The first
  small Chv{\'a}tal closure of $\Frac(G)$, $S_G^{(1)}$, is determined
  by the non-negativity constraints, edge constraints and the rank
  inequalities $\xx(H) \leq \alpha(H)$ as $H = (V_H,E_H)$ varies over
  all induced subgraphs in $G$ containing an odd circuit with vertex
  set $V_H$.
\end{corollary}

It is known that $Q_G^{(1)}$, the first Chv{\'a}tal closure of
$\Frac(G)$, is cut out by the non-negativity, edge and odd circuit
constraints \cite[p. 1099]{SchB}. If an odd circuit is not a hole,
then the corresponding constraint is redundant even though it is
tight. Keeping all odd circuit constraints in $Q_G^{(1)}$, by
Proposition~\ref{prop:A1}, $\A_G^{(1)}$ is precisely the set of
normals of the inequalities describing both $Q_G^{(1)}$ and
$S_G^{(1)}$. However, by Corollary~\ref{cor:oddcircuit normal}, the
right-hand-sides may differ, and $S_G^{(1)}$ could be strictly
contained in $Q_G^{(1)}$.

\begin{example} \label{ex:Kn} Let $G = K_5$. Then $Q_{K_5}^{(1)}$ is
  cut out by the inequalities of $\Frac(K_5)$ along with the 10
  circuit inequalities from the triangles in $K_5$. Its six fractional
  vertices are:
  $$(0,\frac{1}{3},\frac{1}{3},\frac{1}{3},\frac{1}{3}),
  (\frac{1}{3},0,\frac{1}{3},\frac{1}{3},\frac{1}{3}), \ldots, 
  (\frac{1}{3},\frac{1}{3},\frac{1}{3},\frac{1}{3},0), 
  (\frac{1}{3},\frac{1}{3},\frac{1}{3},\frac{1}{3},\frac{1}{3}).$$ 
On the other hand, $S_{K_5}^{(1)}$ is cut out by all the inequalities
  describing $Q_{K_5}^{(1)}$ along with the clique
  inequality $x_1+x_2+x_3+x_4+x_5 \leq 1$, making $S_{K_5}^{(1)}$ equal to
  $\stab(K_5)$.

  In fact, Chv{\'a}tal has shown that the Chv{\'a}tal rank of
  $\Frac(K_n)$ is about $\textup{log} \,n$ \cite{Chv73}. By
  Corollary~\ref{cor:oddcircuit normal}, if $n \geq 3$ is odd, then
  $\scr(\Frac(K_n)) = 1$ since $\stab(K_n)$ is described by the
  inequalities of $\Frac(K_n)$ along with the $n$-clique inequality
  $\sum_{i=1}^{n} x_i \leq 1$.
\end{example}

\begin{corollary} \label{cor:scr=0} For a graph $G$, let $A_G$ denote
  the matrix whose rows are the elements of $\A_G$. Then the following
  are equivalent:
\begin{enumerate}
\item $G$ is bipartite
\item $\Frac(G) = \stab(G)$ (Chv{\'a}tal rank of $\Frac(G)$ is zero) 
\item $\scr(\Frac(G)) = 0$
\item $\scr(A_G) = 0$.
\end{enumerate}
\end{corollary}

\begin{proof} For (1) $\Leftrightarrow$ (2) recall that $G$ is
  bipartite if and only if $G$ has no odd circuits, which is equivalent
  to $\Frac(G) = Q_G^{(1)} = \stab(G)$. Using Lemma~\ref{lem:poly
    characterization of scr} and the fact that $\Frac(G)$ is tight, we
  get (2) $\Leftrightarrow$ (3).  Since $\Frac(G)$ is one polyhedron
  of the form $\{\xx \,:\, A_G \xx \leq \bb \}$, (4) $\Rightarrow$
  (3). On the other hand, Proposition~\ref{prop:A1} shows that if $G$
  is bipartite, then $\A_G^{(1)} = \A_G^{(0)}$ which means $\scr(A_G)
  = 0$ and so (1) $\Rightarrow$ (4).
 \end{proof}

Note that (3) $\Leftrightarrow$ (4) in Corollary~\ref{cor:scr=0} is
highly unusual for a matrix $A$.

\begin{definition} A graph $G$ is {\bf $t$-perfect} if $\stab(G) =
  Q_G^{(1)}$.
\end{definition}

By definition, $t$-perfect graphs are those graphs for which the
Chv{\'a}tal rank of $\Frac(G)$ is one. These graphs have many special
properties and admit a polynomial time algorithm for the stable set
problem. However, no graph theoretic characterization of $t$-perfect
graphs is known. See \cite[Chapter 68]{SchB} for more
details. Example~\ref{ex:Kn} shows that the set of graphs for which
$\scr(\Frac(G)) = 1$ is strictly larger than the set of $t$-perfect
graphs, which raises the following question.

\begin{problem} Characterize the graphs $G$ for which $S_G^{(1)} =
  \stab(G)$, or equivalently, $\scr(\Frac(G)) = 1$.
\end{problem}

We now examine the structure of the vectors in $\A_G^{(2)}$. By a {\em
  cycle} in a graph we mean a collection of circuits in the graph.

\begin{theorem} \label{thm:A2} Every basis $\B \subset \A_G^{(1)}$
  that contributes a vector $\vv$ to $\A_G^{(2)} \setminus \A_G^{(1)}$
  has associated with it a cycle in the hypergraph $G' := ([n],E')$
  where $E'$ is the collection of edges and odd circuits in $G$.  (Two
  hyperedges are adjacent if they share a vertex.)
\end{theorem}

\begin{proof}
  Let $\vv \in \A_G^{(2)} \setminus \A_G^{(1)}$ and $\B = \{\bb_1,
  \ldots, \bb_n\}$ be a basis in $\A_G^{(1)}$ such that $\vv \in
  \Hilb(\B)$. Then there is a $\lambda \in [0,1)^n$ such that $\vv =
  \lambda B$ where $B$ is the $n \times n$ matrix with rows $\bb_1,
  \ldots, \bb_n$. If $p$ of the elements in $\B$ are $-\ee_i$'s, then
  $p < n$ and we may assume that $\bb_i = -\ee_i$ for $i= n-p+1,
  \ldots, n$. We will show that $\B$, and hence $\vv$, can be
  associated with a cycle in $G'$.

  Let $B'$ be the top left $(n-p) \times (n-p)$ submatrix of $B$. 
  Then by Proposition~\ref{prop:A1},  
  $\bb_1, \ldots, \bb_{n-p}$ are all
  characteristic vectors of edges and odd circuits in $G$, and hence
  $B' \in \{0,1\}^{(n-p) \times (n-p)}$. Consider the $j$-th column
    in $B'$. This column is not all zero since $\B$ is a basis. If it
    has exactly one $+1$, then $v_j = 0 = \lambda_j$, and we may
    ignore the $j$-th row and column of $B$. Therefore, assume that
    each column of $B'$ has at least two $+1$'s. Each row of $B'$ has
    at least one $+1$, since otherwise, $\textup{det}(B') = 0 =
    \textup{det}(B)$. Suppose there are $q$ rows in $B'$ with exactly
    one $+1$. By permuting rows and columns in $B'$, we may assume
    that these rows are at the bottom of $B'$ and that they contribute
    a $q \times q$ identity matrix in the bottom right of $B'$.  If $n
    = q+p$ then $|\textup{det}(B)| = 1$ and there is no $\vv$ as above
    to consider. Therefore, $n > q+p$. (The structure of $B$ is shown
    below where $\delta$ is used for an entry that may be $0$ or $1$.)

    Let $B''$ denote the top left $(n-q-p) \times (n-q-p)$ submatrix
    of $B'$. By the same argument as for $B'$, each column of $B''$
    has at least two $+1$'s. Counting the $+1$'s in $B''$, each row of
    $B''$ must also have at least two $+1$'s. Let the vertices
    indexing the columns of $B''$ be $V''$. Then each $v \in V''$ is
    incident to at least two hyperedges in $G'$ from the set of
    hyperedges indexed by $\bb_1, \ldots, \bb_{n-q-p}$. This implies
    that there exists a circuit or collection of circuits in $G'$
    through the vertices in $V''$ using the above hyperedges.
\begin{tiny}
\begin{displaymath}
B = \left( \begin{array}{c c c c c | c c c c | c c c c}
\mathbf{1} & \mathbf{1} & \mathbf{\delta} & & \mathbf{\delta} & \delta & \delta & & \delta & \delta & \delta & & \delta \\
\mathbf{\delta} & \mathbf{1} & \mathbf{1} & & \mathbf{\delta} & \delta & \delta & & \delta & \delta & \delta & & \delta \\
& & & \ddots & & & & & & & & & \\
\mathbf{1} & \mathbf{\delta} & \mathbf{\delta} & & \mathbf{1} & \delta & \delta & & \delta & \delta & \delta && \delta \\
\hline
0 & 0 & 0 & & 0 & \mathbf{1} & \mathbf{0} & & \mathbf{0} & \delta & \delta & & \delta \\
0 & 0 & 0 & & 0 & \mathbf{0} & \mathbf{1} & & \mathbf{0} & \delta & \delta &  & \delta \\
& & & & & & & \ddots & & & & & \\
0 & 0 & 0 & & 0 & \mathbf{0} & \mathbf{0} & & \mathbf{1} & \delta & \delta & & \delta \\
\hline
0 & 0 & 0 & & 0 & 0 & 0 & & 0 & \mathbf{-1} & \mathbf{0} & & \mathbf{0}\\
0 & 0 & 0 & & 0 & 0 & 0 & & 0 & \mathbf{0} & \mathbf{-1} & & \mathbf{0}\\
& & & & & & & & & & & \ddots & \\
0 & 0 & 0 & & 0 & 0 & 0 & & 0 & \mathbf{0} &\mathbf{0} & & \mathbf{-1}
\end{array} \right)
\end{displaymath}
\end{tiny}
 \end{proof}

In the rest of this section we will show that the normals of many
complicated families of valid inequalities for $\stab(G)$ appear in
$\A_G^{(2)}$ which shows that they are all derived from hypercycles in
$G'$ as in Theorem~\ref{thm:A2}. 

\begin{example} \label{ex:clawfree} ({\bf Claw-free graphs}) A graph
  $G$ is claw-free if it does not contain an induced $K_{1,3}$ ({\em
    claw}). It is known that the maximum stable set problem in a
  claw-free graph can be solved in strongly polynomial time
  \cite[Chapter 69]{SchB}, but it is a long-standing open problem to
  give a description of $\stab(G)$. Claw-free graphs have been shown
  to have complicated facet inequalities \cite{Giles-Trotter}, \cite{LOSS},
  and a full characterization of their rank facet inequalities is also
  known \cite{GalluccioSassano}.
 
  It was shown in \cite{Giles-Trotter} that for a fixed positive
  integer $a$, there are claw-free graphs on $n := 2a(a+2)+1$ vertices
  that have a facet normal with coefficients $a$ and $a+1$. The
  corresponding facet inequalities are produced in one round of the
  Chv{\'a}tal procedure if one starts with the clique and
  non-negativity constraints. Therefore, these normals appear in at
  most three rounds of IBN, since clique normals appear in two. It is
  not hard to see that these normals are produced in two rounds of
  IBN.

To illustrate Theorem~\ref{thm:A2}, we pick the example on pp. 321 of
\cite{Giles-Trotter} which considers the claw-free graph $G$ that is
the complement of the graph in Figure~\ref{fig:fig4}.  In this case,
$\stab(G)$ has 35 facets and the following is an example of a facet
inequality with more than two non-zero coefficients.
$$ 2x_1+2x_2+2x_3+2x_4+2x_5+x6+x_7+3x_8+x_9+3x_{10} \leq 4 $$
After permuting coordinates to be in the order
$[6,7,9,1,2,3,4,5,8,10]$, the normal of the above inequality is $\vv =
(1,1,1,2,2,2,2,2,3,3)$ and it lies in $\Hilb(\B)$, where $\B$ is a
basis in $\A_G^{(1)}$ for which $B$ is as follows: 
\begin{displaymath}
B = \left( \begin{array}{rrr|rrrrr|rr}
{\bf 1} & {\bf 1} & 0 & 1 & 1 & 1 & 1 & 1 & 1 & 1 \\
{\bf 1} & 0 & {\bf 1} & 1 & 1 & 1 & 1 & 1 & 1 & 1 \\
0 & {\bf 1} & {\bf 1} & 1 & 1 & 1 & 1 & 1 & 1 & 1 \\
\hline
0 & 0 & 0 & 1 & 0 & 0 & 0 & 0 & 0 & 1 \\
0 & 0 & 0 & 0 & 1 & 0 & 0 & 0 & 0 & 1 \\
0 & 0 & 0 & 0 & 0 & 1 & 0 & 0 & 1 & 0 \\
0 & 0 & 0 & 0 & 0 & 0 & 1 & 0 & 1 & 0 \\
0 & 0 & 0 & 0 & 0 & 0 & 0 & 0 & 1 & 1 \\
\hline
0 & 0 & 0 & 0 & 0 & 0 & 0 & 0 & -1 & 0 \\
0 & 0 & 0 & 0 & 0 & 0 & 0 & 0 & 0 & -1
\end{array} \right).
\end{displaymath}
Check that $\vv =
(\frac{1}{2},\frac{1}{2},\frac{1}{2},\frac{1}{2},\frac{1}{2},
\frac{1}{2},\frac{1}{2},\frac{1}{2},0,0)B$ and notice the
hypertriangle through the vertices $6,7,9$ made up of 
three 9-circuits in $G$.

\begin{figure} 
\label{fig:fig4}
\includegraphics[width=3in]{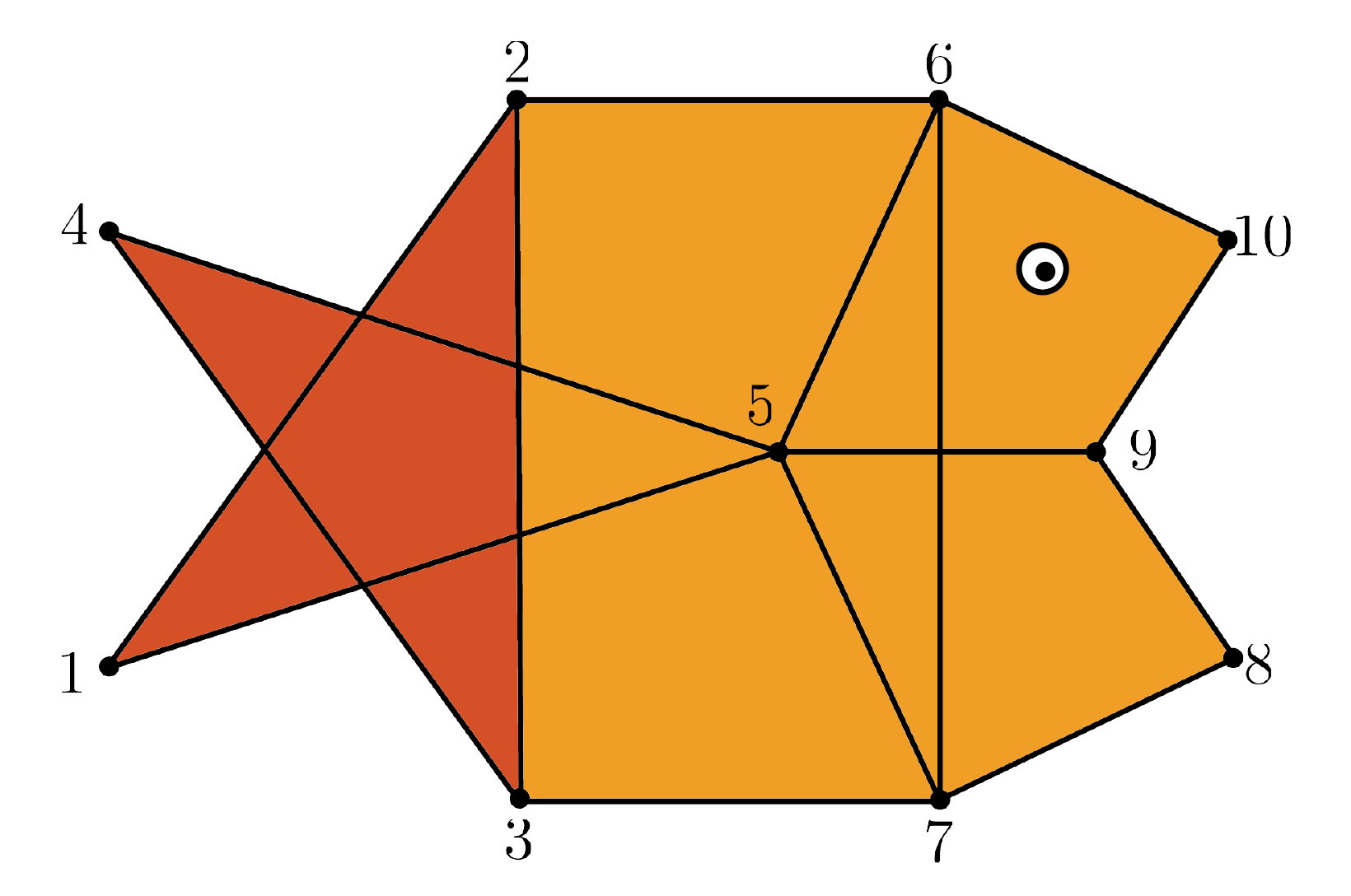}
\caption{Fig. 4 in \cite{Giles-Trotter}}
\end{figure}

In \cite{LOSS}, the authors extend the above example from
\cite{Giles-Trotter} by showing two claw-free graphs called ``fish in
a net'' and ``fish in a net with bubble'' each with a facet inequality
that has several different non-zero coefficients. Both normals appear
in $\A_G^{(2)}$. We illustrate the ``fish in a net with bubble''
case. Let $G$ be the complement of the graph shown in
Figure~\ref{fig:fish}.  Then $\stab(G)$ has the following facet
inequality \cite{LOSS}:
\begin{eqnarray*}
3x_1+3x_2+3x_3+3x_4+4x_5+4x_6+4x_7+5x_8+4x_9 &\\
+5x_{10}+5x_{11}+4x_{12}+6x_{13}+2x_{14}+2x_{15}+2x_{16}+6x_{17} & \leq 8
\end{eqnarray*}

\begin{figure} 
\includegraphics[width=\textwidth]{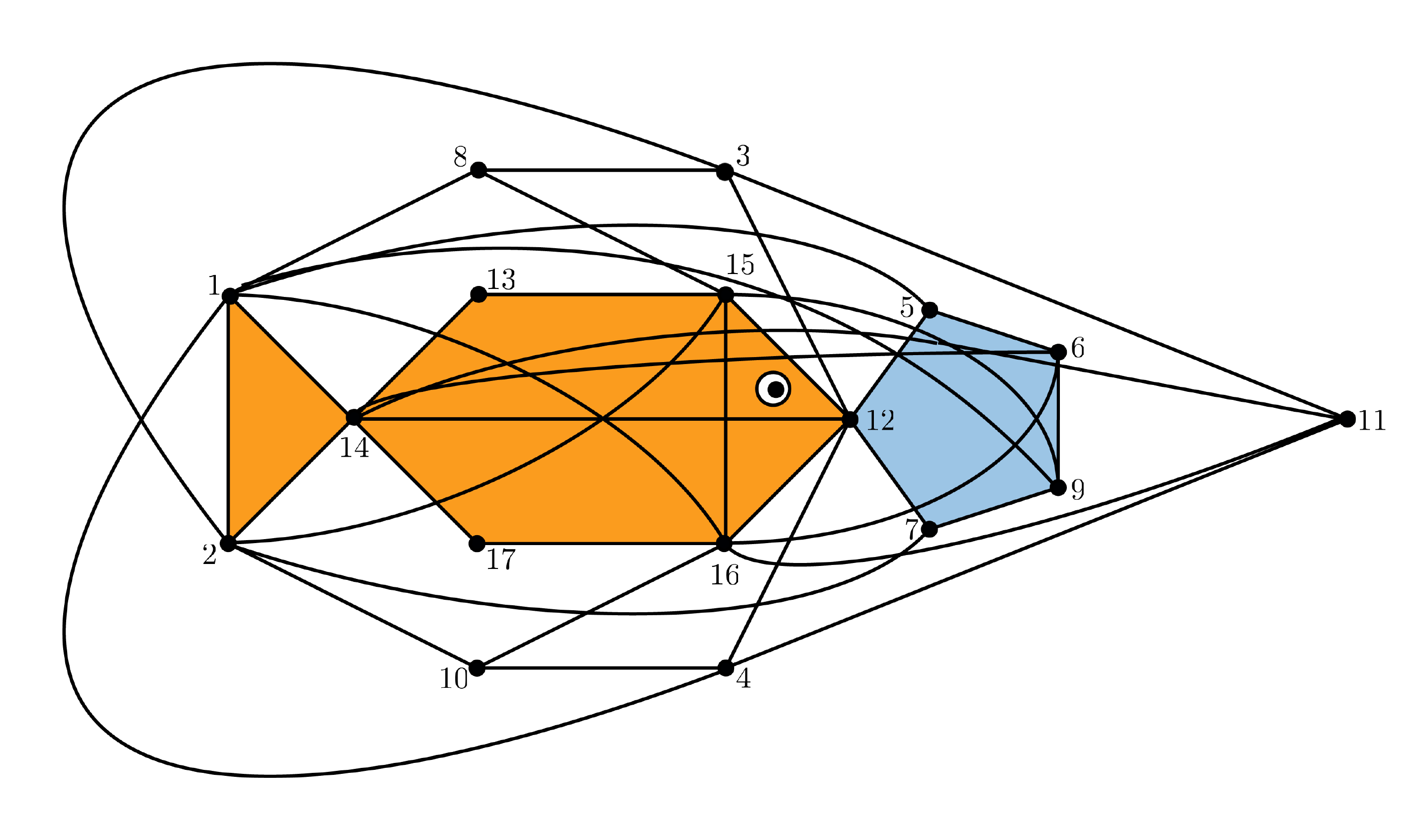}
\caption{Fish in a net with bubble}
\label{fig:fish}
\end{figure}

Let $\B$ be the rows of the following $17 \times 17$ matrix $B$. 
\begin{displaymath}
\left( \begin{array}{c c c c c c c c c c c c c | r r r r}
\mathbf{1} & 1 & 1 & 1 & 0 & 1 & 1 & 1 & 1 & 1 & 1 & 0 & \mathbf{1} & 1 & 1 & 1 & 1\\
\mathbf{1} & \mathbf{1} & 1 & 1 & 1 & 0 & 1 & 1 & 0 & 1 & 1 & 1 & 1 & 1 & 1 & 1 & 1\\
1 & \mathbf{1} & \mathbf{1} & 1 & 1 & 1 & 0 & 1 & 1 & 1 & 1 & 0 & 1 & 1 & 1 & 1 & 1\\
0 & 1 & \mathbf{1} & \mathbf{1} & 1 & 1 & 1 & 1 & 1 & 1 & 1 & 1 & 1 & 0 & 0 & 0 & 1\\
1 & 0 & 1 & \mathbf{1} & \mathbf{1} & 1 & 1 & 1 & 1 & 1 & 1 & 1 & 1 & 0 & 0 & 0 & 1\\
1 & 1 & 0 & 1 & \mathbf{1} & \mathbf{1} & 1 & 1 & 1 & 1 & 1 & 1 & 1 & 0 & 0 & 0 & 1\\
1 & 1 & 1 & 0 & 1 & \mathbf{1} & \mathbf{1} & 1 & 1 & 1 & 1 & 1 & 1 & 0 & 0 & 0 & 1\\
0 & 0 & 0 & 0 & 1 & 1 & \mathbf{1} & \mathbf{1} & 0 & 1 & 1 & 1 & 1 & 0 & 0 & 0 & 1\\
0 & 0 & 0 & 0 & 0 & 0 & 0 & \mathbf{1} & \mathbf{1} & 0 & 0 & 1 & 1 & 0 & 0 & 0 & 1\\
0 & 0 & 0 & 0 & 0 & 0 & 0 & 0 & \mathbf{1} & \mathbf{1} & 0 & 0 & 0 & 0 & 0 & 0 & 1\\
0 & 0 & 0 & 0 & 1 & 1 & 1 & 1 & 1 & \mathbf{1} &\mathbf{1} & 0 & 1 & 0 & 0 & 0 & 1\\
0 & 0 & 0 & 0 & 0 & 0 & 0 & 0 & 0 & 0 & \mathbf{1} & \mathbf{1} & 1 & 0 & 0 & 0 & 0\\
0 & 0 & 0 & 0 & 0 & 0 & 0 & 0 & 0 & 0 & 0 & \mathbf{1} & \mathbf{1} & 0 & 0 & 0 & 1\\
\hline
0 & 0 & 0 & 0 & 0 & 0 & 0 & 0 & 0 & 0 & 0 & 0 & 0 & \mathbf{-1} & 0 & 0 & 0\\
0 & 0 & 0 & 0 & 0 & 0 & 0 & 0 & 0 & 0 & 0 & 0 & 0 & 0 & \mathbf{-1} & 0 & 0\\
0 & 0 & 0 & 0 & 0 & 0 & 0 & 0 & 0 & 0 & 0 & 0 & 0 & 0 & 0 & \mathbf{-1} & 0\\
0 & 0 & 0 & 0 & 0 & 0 & 0 & 0 & 0 & 0 & 0 & 0 & 0 & 0 & 0 & 0 & \mathbf{-1}\\
\end{array} \right)
\end{displaymath}

Then $\bb_1, \ldots, \bb_{13}$ are characteristic vectors of odd
circuits in $G$. We denote the consecutive vertices in one such odd
circuit for each $\bb_i$ in the table below.

\begin{displaymath}
\begin{array}{ll}
\bb_1 & 1 , 3 , 4 , 2 , 6 , 7 , 8 , 9 , 10 , 11 , 13 , 16
, 14 , 15 , 17 \\ 
\bb_2 & 1 , 3 , 2 , 5 , 4
, 7 , 8 , 10 , 11 , 12 , 13 , 16 , 14 , 15
, 17\\ 
\bb_3 & 1 , 3 , 2 , 5 , 4 , 6 , 8 , 9 , 10 
, 11 , 15 , 14 , 16 , 13
, 17 \\
\bb_4 & 2 , 3 , 4 , 5 , 7 , 6 , 8 , 9 , 10 
, 11 , 12 , 13 , 17\\ 
\bb_5 & 1 , 3 , 4 , 5 , 7 , 6 , 8 , 9 , 10 
, 11 , 12 , 13 , 17\\
\bb_6 & 1 , 4 , 2 , 5 , 7 , 6 , 8 , 9 , 10 
, 11 , 12 , 13 , 17\\
\bb_7 & 1 , 3 , 2 , 5 , 7
, 6 , 8 , 9 , 10 
, 11 , 12 , 13 , 17 \\ 
\bb_8 & 5 , 7 , 6 , 8 , 10
, 11 , 12 , 13 , 17\\
\bb_9 & 8 , 9 , 12 , 13 , 17 \\
\bb_{10} & 9 , 10 , 17\\ 
\bb_{11} & 5 , 7 , 6 , 8 , 9
, 10 , 11 , 13 , 17 \\
\bb_{12} & 11 , 12 , 13 , 11\\
\bb_{13} & 12 , 13 , 17 , 12 
\end{array}
\end{displaymath}

Therefore, $\B \subset \A_G^{(1)}$ and check that $\textup{det}(B) =
18$, hence, $\B$ is a basis. It can be verified using a Hilbert basis
package such as Normaliz \cite{Normaliz} that the normal {\bf v} from
the facet inequality above is in $\textup{Hilb}(\B)$: $$\vv =
\left( \frac{2}{3}, \frac{2}{3}, \frac{2}{3}, \frac{1}{3}, \frac{1}{3},
\frac{1}{3}, \frac{1}{3}, \frac{2}{3}, \frac{1}{3}, \frac{1}{3},
\frac{2}{3}, \frac{1}{3}, \frac{2}{3}, 0, 0, 0, 0 \right) B,$$ which proves
that $\vv \in \A_G^{(2)}$. The hypercycle associated with $\B$ is
indicated by the bold $1$'s in the matrix $B$.
\end{example}

It is a long-standing open problem to give a complete linear
inequality description of $\stab(G)$ when $G$ is a claw-free
graph. The following would be a step toward settling this problem.

\begin{problem} Is $\scr(\Frac(G)) \leq 2$ for all claw-free graphs $G$?
\end{problem}

We now derive various corollaries to Theorem~\ref{thm:A2}.

\begin{corollary}
  The normals of all clique inequalities lie in $\A_G^{(2)}$.
\end{corollary}

\begin{proof} Corollary~\ref{cor:oddcircuit normal} showed that the
  normals of all odd clique inequalities lie in $\A_G^{(1)}$. Suppose
  $H$ is an even clique in $G$ with vertex set $V_H$. For $i,j \in
  V_H$, let $C_i$ be an odd circuit through all vertices of $V_H$
  except $i$ and similarly, $C_j$ be an odd circuit through all
  vertices of $V_H$ except $j$. Then $\ee(C_i), \ee(C_j)$ and
  $\ee_i+\ee_j$ are all present in $\A_G^{(1)}$ by
  Proposition~\ref{prop:A1}. The odd circuits $C_i,C_j$ and the edge
  $\{i,j\}$ together form a triangle in the hypergraph $G'$, and the
  vectors $\ee(C_i), \ee(C_j)$ and $\ee_i+\ee_j$ are linearly
  independent since for any $k \in V_H \setminus \{i,j\}$, the $3
  \times 3$ submatrix indexed by $i,j,k$, of the $3 \times n$ matrix
  whose rows are these three vectors is non-singular. Dividing 
  the sum of the three vectors by $2$ produces
  $\ee(H)$. This vector is in the minimal Hilbert basis of the cone
  spanned by $\ee(C_i), \ee(C_j)$ and $\ee_i+\ee_j$ since its
  restriction $(1,1,1)$ to the coordinates indexed by $i,j,k$ is in
  the minimal Hilbert basis of the cone spanned by the same
  restriction of $\ee(C_i), \ee(C_j)$ and $\ee_i+\ee_j$. Since
  $\ee(C_i)$, $\ee(C_j)$, $\ee_i + \ee_j$ can be extended to a basis
  in $\A_G^{(1)}$, the result follows.
 \end{proof}

\begin{definition} \label{def:perfect, hperfect}
\begin{enumerate}
\item A graph $G = ([n],E)$ is {\bf perfect} if $\stab(G)$ is cut out
  by the non-negativity and clique inequalities.
\item A graph $G = ([n],E)$ is {\bf h-perfect} if $\stab(G)$ is cut
  out by the non-negativity, odd circuit and clique inequalities.
\end{enumerate}
\end{definition}

All perfect graphs are $h$-perfect. Many well-known
classes of graphs such as bipartite, comparability and chordal graphs
are perfect \cite{SchB}.

\begin{corollary} \label{cor:scr of perfect and h-perfect} If $G$ is
  h-perfect then $\scr(\Frac(G)) = 2$.
\end{corollary}

If $H$ is a subgraph of $G$, then note that the configuration
$\A_H^{(0)}$ is a subset of $\A_G^{(0)}$ after padding all coordinates
corresponding to vertices of $G$ that are not in $H$ by zeros. This implies
that $\A_H^{(k)}$ is also a subset of $\A_G^{(k)}$, for any positive
integer $k$, after the same padding by zeros. Therefore, if we need to
show that a facet normal of $\stab(G)$ whose support lies in the
vertices of $H$ appears in $\A_G^{(k)}$, then it suffices to show that
it appears in $\A_H^{(k)}$. 

\begin{corollary} Normals of antihole and odd wheel inequalities 
  appear in $\A_G^{(2)}$.
\end{corollary}

\begin{proof}
  By the above discussion, we may assume without loss of generality
  that $G$ is an antihole. Since the vertices of an odd antihole
  support an odd circuit, its characteristic vector appears in
  $\A_G^{(1)}$. If $G$ is an even antihole, then $G$ contains $n$ odd
  circuits each going through all vertices of $G$ except one. The $n
  \times n$ matrix $B$ whose rows are the characteristic vectors of
  these odd circuits has all diagonal entries equal to zero and all
  off-diagonal entries equal to one. Since $B$ is non-singular,
  its rows form a basis in $\A_G^{(1)}$. Dividing the sum of the rows
  of $B$ by $n-1$ produces $\ee(G)$, which is the
  unique new element in the minimal Hilbert basis of
  $\cone(\B)$. 

  Again assume without loss of generality that $G$ is an odd wheel
  with central vertex $u_0$ and remaining vertices $u_1, \cdots,
  u_{2k-1}$. Let $B$ be the $2k \times 2k$ matrix whose rows are the
  characteristic vectors of the $2k-1$ triangles in $G$ and
  $-\ee_0$. Then $B$ is non-singular and that half the sum
  of its rows is the normal of the odd wheel inequality.

\begin{displaymath}
 \left( \begin{array}{r r r r r r r }
1  & 1 & 1 & 0 & 0 & \ldots & 0 \\
1  & 0 & 1 & 1 & 0 & \ldots & 0\\
1  & 0 & 0 & 1 & 1 & \ldots & 0\\
  &   &     &    & & \ddots&  \\
1  & 1 & 0 & 0 & 0 & \ldots & 1\\ 
-1 & 0 & 0 & 0 & 0 & \ldots & 0
 \end{array} \right)
 \end{displaymath}

The $2k-1$ triangles are the hyperedges that form an odd cycle in $G'$
which underlies this normal.
 \end{proof}

Recall that the {\bf line graph}, $L(G)$, of a graph $G=(V,E)$ is the
graph $L(G) = (E,F)$ where $\{e,e'\} \in F$, for $e,e' \in E$, if and
only if $e$ and $e'$ share a vertex in $G$.  A complete linear
description of $\stab(L(G))$ was given by Edmonds as follows (see
\cite[p. 440]{SchA}):

$$\stab(L(G)) = \left \{ \xx \in \RR^E \,:\, 
\begin{array}{ll}
x_e \geq 0 & \forall \,\,e \in E \\
\sum_{v \in e} x_e \leq 1 & \forall \,\, v \in V\\
\sum_{e \in E[U]} x_e  \leq \lfloor \frac{|U|}{2} \rfloor & \forall
\,\,U \subseteq V, \, |U| \,\,\textup{odd} 
\end{array}
\right \},$$
where $E[U]$ denotes the edges in $E$ that have both end points in
$U$. Note that the second class of inequalities in the
description of $\stab(L(G))$ are clique inequalities.  

\begin{corollary} \label{cor:line graphs} For any graph $G$,
  $\scr(\Frac(L(G))) \leq 3$.
\end{corollary}

\begin{proof}
  It is known that $\stab(L(G))$ is the first Chv{\'a}tal closure of
  the polytope described by the clique and non-negativity constraints
  from $L(G)$. Since clique normals are in $\A_{L(G)}^{(2)}$, it
  follows that $\scr(\Frac(L(G))) \leq 3$.
 \end{proof}

\section{Lower bounds} \label {sec:lowerbounds} 

In this section we establish lower bounds on SCR in various
situations. We also discuss computational evidence that supports
possible upper bounds in some of these cases. Note that by
Proposition~\ref{prop:ranks}, a lower bound on the SCR or Chv{\'a}tal
rank of a system $A \xx \leq \bb$ is also a lower bound on the
corresponding rank of $A$. On the other hand, an upper bound on either
rank of $A$ is an upper bound on the corresponding rank of $A \xx \leq
\bb$ for any $\bb$.

\begin{theorem}\label{thm:largeSCR} For $m,n \geq 3$, the small 
Chv{\'a}tal rank of $A \xx \leq \bb$ (and hence of $A$) can grow 
exponentially in the size of the input.
\end{theorem}

In proving Theorem~\ref{thm:largeSCR}, we may assume 
$m=n=3$. All other cases follow by
adjoining inequalities that do not affect Chv{\'a}tal rank or SCR. 
Let $j \geq 2$ be arbitrary and set
$$ A = \left( \begin{array}{ccc}
    1 & 0 & 0 \\
    0 & 1 & 0 \\
    1 & j & 2j-1 \end{array} \right). $$
We will show that the SCR of $A$ is $j-1$ which is exponential in the 
bit size of $A$.  To do this, we explicitly describe $\A^{(k)}$ 
for all $k$ and prove that $\A^{(j-1)} = \A^{(j)}$,
so the SCR of $A$ is at most $j-1$. Then we identify a 
vector in $\A^{(j-1)} \setminus \A^{(j-2)}$ that is a facet normal
of an integer hull $Q_{\bb}^I$. 
\vspace{0.2cm}
 
For $1 \, \leq \, k \, \leq \, j-1$, 
define an integral polygon 
\begin{eqnarray*} R^k & := & \conv \{ (0,0), (k+1,k+1),(j,2j-1-k), (j,
  2j-1) \} \subseteq \RR^2_{\geq 0} \\
& = & \{ (x,y) \in \RR^2: \; x \leq y,\; 2x \leq y + k + 1,\;
x  \leq  j,\; y \leq \left( \frac{2j-1}{j} \right) x \}. \end{eqnarray*}
For $k = j-1$, the second and third points in the convex hull
description coincide; for $k < j-1$
the four points are distinct and in convex position.  

\begin{figure}
\label{fig:A_65} 
\begin{center}
 \includegraphics[height=5.0cm]{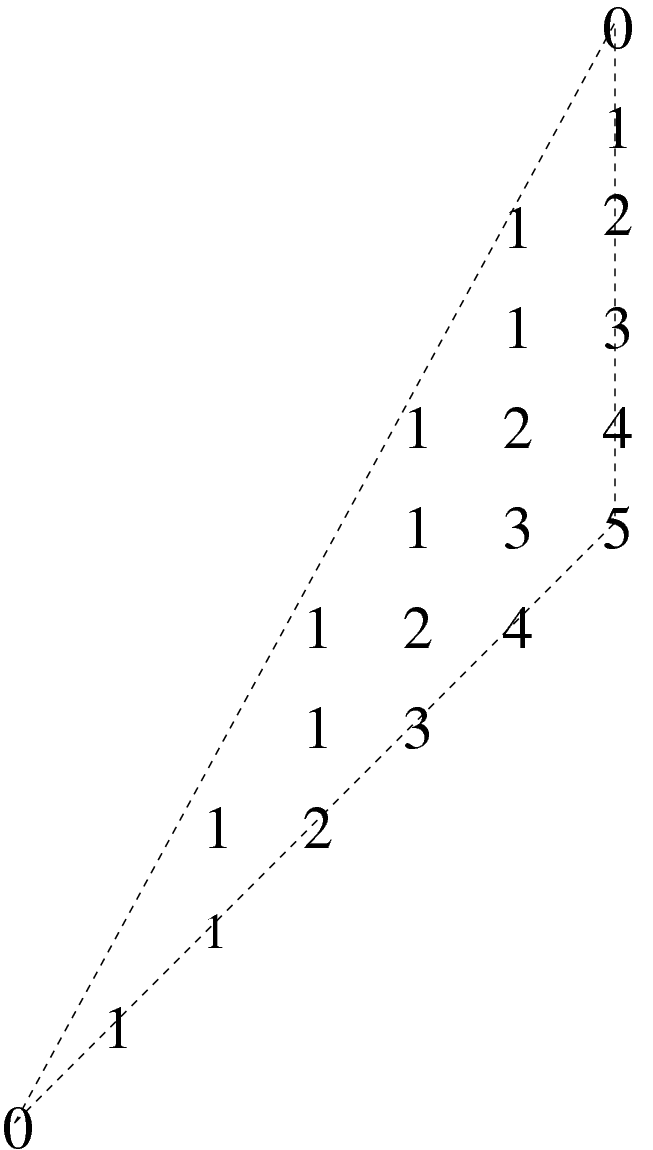}
 \end{center}
 \caption{The polygon $R^{j-1}$ (with $j=6$) used
 to prove Theorem~\ref{thm:largeSCR}.  Each integer point $(a,b)$ 
 in the polygon is labeled by the smallest $k$ such that 
 $(1,a,b)$ appears in $\A_6^{(k)}$.}
 \end{figure}

\begin{lemma} \cite[Proposition 5.1]{HMS} \label{lem:polygonsupernormal} 
  Let $R$ be an integral polygon in $\RR^2$.  The configuration 
   $\{ (1,a,b): \, (a,b) \in R \cap \ZZ^2 \}$ is supernormal.
\end{lemma}

\begin{lemma} \label{lem:Ajk} For $1 \leq k \leq j-1$, 
$$\A^{(k)} = \{(0,1,0)\} \cup \{ (1,a,b): \, (a,b) \in R^k \cap
\ZZ^2\}.$$ 
\end{lemma}

\begin{proof} Induct on $k$. For $k=1$, we have 
  $$R^1 = \{ (x,y) \in \RR^2: \; x \leq y,\; 2x \leq y + 2,\;
  x  \leq  j,\; y \leq \left( \frac{2j-1}{j} \right) x \} $$
and it is easy to check that 
$$R^1 \cap \ZZ^2= \{(0,0)\} \cup \{(i,2i-1) \,:\, 1 \leq i \leq j\} \cup
  \{(i,2i-2) \,:\, 2 \leq i \leq j\}.$$

Observe that
  $$
  (1,i,2i-1) = \left(\frac{2j-2i}{2j-1},\frac{j-i}{2j-1},
  \frac{2i-1}{2j-1}\right) A$$
  for $1 \leq i \leq j$ and that
  $$
  (1,i,2i-2) = \left(\frac{2j-2i+1}{2j-1},\frac{2j-i}{2j-1},
  \frac{2i-2}{2j-1}\right) A$$
  for $2 \leq i \leq j$, so all the points
  in $\{1\} \times (R^1 \cap \ZZ^2)$ are in the fundamental
  parallelepiped of $\A$. Since all the first coordinates are one,
  no element of $\{1\} \times R^1 \cap \ZZ^2$ is a sum of others.
  Also, no two elements of $\{1\} \times R^1 \cap \ZZ^2$ differ by a
  multiple of $(0,1,0)$.  Thus $\{1\} \times (R^1
  \cap \ZZ^2) \subseteq \Hilb(\A) \subseteq \A^{(1)}.$ 

On the other hand, if $h = c_1 (1,0,0) + c_2 (0,1,0) + c_3
  (1,j,2j-1)$ is an integer point in the fundamental parallelepiped of
  $\A$ (so $0 \leq c_1, c_2, c_3 < 1$), then $c_3 = \frac{p}{2j-1}$
  for some integer $1 \leq p \leq 2j-2$ and $c_1$ and $c_2$ are
  uniquely determined by $c_3$, so $h$ must be one of the listed 
  points in $R^1 \cap \ZZ^2$. 

For the induction step, first assume that $\A^{(k-1)}$ 
  contains $\{1\} \times R^{k-1} \cap \ZZ^2$ for some $k \geq 2$.  The 
  difference between $R^{k-1}$ and $R^k$ is that the
  inequality $2x \leq y+k$ is relaxed to $2x \leq y+k+1$. So we must 
show that the new vectors in
  $\A^{(k)}$ include
\begin{equation} \label{eq:newinAjk} 
\{ (1,k+i,k+2i-1): 1 \, \leq i \, \leq j-k \}.
\end{equation}   
For each $1 \, \leq \, i \, \leq \, j-k$, the three vectors $(0,1,0)$,
$(1, k+i-1, k+2i-2)$, and $(1,k+i, k+2i)$ appear in $\A^{(k-1)}$ by
the induction hypothesis.  The basis cone $C$ that they span has normalized
volume two, and $(1,k+i,k+2i-1)$ (half the sum of the three vectors) 
is the unique integer point in the interior of the fundamental 
parallelepiped.  Thus $(1,k+i,k+2i-1) \in \Hilb(C) \subseteq \A^{(k)}.$ 

Next assume for some $k$ that $\A^{(k-1)}$ contains no other vectors. 
By Lemma~\ref{lem:polygonsupernormal}, the previous paragraph, and the 
induction hypothesis, the set 
$$\A_{j}^{(k-1)} \setminus \{ (0,1,0) \} \; = \; R^k \cap \ZZ^2$$ 
is supernormal.  Thus the only bases of $\A^{(k-1)}$ that might 
contribute new vectors to $\A^{(k)}$ are those that include $(0,1,0).$ 
Any new vector obtained this way would be of the form 
$(1,a,b)$ for $(a-1,b)$ strictly in the interior of 
$R^{k-1}$ and $(a,b)$ outside $R^{k-1}$.
From the inequality description of $R^{k-1}$, this vector must indeed
be of the form~(\ref{eq:newinAjk}); see Figure~\ref{fig:A_65}.  
 \end{proof}

\begin{lemma} \label{lem:Ajj-1supernormal}
The configuration $\A^{(j-1)}$ is supernormal.
\end{lemma}

\begin{proof}
  By the same argument as above, any vector $\vv \in \A^{(j)} \setminus 
\A^{(j-1)}$ is
of the form $(1,a,b)$ for $(a-1, b)$ in the interior of $R^{j-1}$ and
  $(a,b)$ outside $R^{j-1}$. However, $R^{j-1}$ is a triangle whose 
  right boundary  consists only of segments of the line $y = x$ and of
  the line $x \leq j$, so no such $(a,b)$ exists. Thus 
  $\A^{(j)} = \A^{(j-1)}$.
 \end{proof}

\smallskip\noindent {\it Proof of Theorem~\ref{thm:largeSCR}}:
By Lemma~\ref{lem:Ajk}, we have 
$(1,j,j)^t \in \A^{(j-1)} \setminus \A^{(j-2)}$.  So it will suffice
to show that the inequality 
\begin{equation} \label{eq:Ajkhull}  (1,j,j) \, \xx \, \leq \, {\bf 0}
\end{equation}
defines a facet of the integer hull 
$$P_j := \{\xx \in \RR^3 \, : \, A \xx \, \leq \, (0,0,j-1)^t \}^I. $$

Let $\yy = (y_1,y_2,y_3) \in Q_{(0,0,j-1)^t} \cap \ZZ^3$. 
We first show that $\yy$ satisfies
~(\ref{eq:Ajkhull}). If $y_3 \leq 0$, then since we already know $y_1,
y_2 \leq 0$, immediately $\yy$ satisfies ~(\ref{eq:Ajkhull}). If 
$y_3 = 1$ and $y_2 \leq -1$, again $\yy$ satisfies
~(\ref{eq:Ajkhull}).
If $y_3 = 1$ and $y_2 = 0$, then to satisfy the last inequality in $A
  \xx \leq (0,0,j-1)^t$, $y_1 \leq -j$ and again $\yy$ satisfies
  ~(\ref{eq:Ajkhull}).

Finally, suppose $y_3 \geq 2$.  Rewrite $x_1 + jx_2 + (2j-1)x_3 \,
\leq \, j-1$ as
\begin{equation} \label{eq:x3rewrite}
x_1 + j x_2 \, \leq \, (j-1) - x_3(2j-1). \end{equation}
Then 
$$ \begin{array}{ccccc} 
(1,j,j) \, \yy & =    & y_1 + jy_2 + jy_3 & \leq & (j-1) - y_3(2j-1)
+ j y_3 \\
               & =    & j + y_3(1-j) - 1 & \leq & j + 2(1-j) - 1 \\
               & =    & 1-j              & <    & 0 \end{array} $$ 
where the first inequality follows from~(\ref{eq:x3rewrite}), 
the second from $y_3 \geq 2$, and the last from $j \geq 2$.  
Thus the inequality~(\ref{eq:Ajkhull}) is valid on all integer points
of $Q_{(0,0,j-1)^t}$ and hence on $P_j$.

To finish the proof we must argue that ~(\ref{eq:Ajkhull}) is a
facet inequality of $P_j$. This follows from the observation that 
the three affinely independent integer points $(0,-1,1)^t$,
$(0,0,0)^t$, and $(-j,0,1)^t$ in $P_j$ satisfy ~(\ref{eq:Ajkhull})
with equality. \endproof 

Many optimization problems are modeled as $0/1$ integer programs, in
which case the starting linear programming relaxation is a polytope in
the unit cube $C_n = [0,1]^n.$ For such polytopes, it is known that
Chv{\'a}tal rank is bounded above by $n^2(1 + \log n)$, and there are
examples with Chv{\'a}tal rank at least $(1+\epsilon) n$
\cite{EisenbrandSchulz}. We will derive a lower bound for SCR of the
same order, using quite different techniques.

\begin{theorem} \label{thm:0/1SCRlowerbound}
There are systems $A \xx \leq \bb$ defining polytopes contained in
the unit cube $C_n$ whose small Chv{\'a}tal ranks are at least $n/2 -
o(n)$.  
\end{theorem}   

\begin{observation} \label{obs:HAnorm} If $\vv \in \Hilb(\{\vv_1,
  \ldots, \vv_n\})$, then $\|\vv\|_\infty < n (\max_i
  \|\vv_i\|_\infty)$ since $\vv = \sum_{i=1}^{n} \lambda_i \vv_i$ for
  $0 \leq \lambda_i < 1$, $i=1,\ldots,n$. 
\end{observation}

\vspace{.2cm} \smallskip\noindent {\it Proof of
  Theorem~\ref{thm:0/1SCRlowerbound}}:
Given any 0/1 polytope $Q$, we
can find a relaxation $P$ contained in $C_n$ and whose facet normals
are 0/1/-1 vectors. For instance, for any $U \subseteq [n]$, 
the inequality 
$$
\sum_{i \in U} x_i - \sum_{i \notin U} x_i \; \leq \; \left| U
\right| -1$$
is violated by $\ee(U)$ but satisfied by every other
vertex of $C_n$. Define $P$ by starting with $C_n$ and
adjoining such an inequality for each vertex of $C_n$ that is not
in $Q$.  

Using a construction by Alon and Vu \cite{AV} of 0/1 matrices with large
determinants, Ziegler \cite[Corollary 26]{Ziegler01} constructs an
$n$-dimensional 0/1 polytope $Q$ with a (relatively prime integer)
facet normal $\vv$ whose $\infty$-norm is at least
$\frac{(n-1)^{(n-1)/2}}{2^{2n + o(n)}}.$ Let $P$ be as above for
this $Q$, and let $k$ be the SCR of the system 
$A \xx \leq \bb$ defining $P$.  By definition, $\vv \in
\A^{(k)}$. Since $\A$ consists entirely of 0/1/-1 vectors, 
we get by repeatedly applying Observation~\ref{obs:HAnorm} that 
$$n^k > \frac{(n-1)^{(n-1)/2}}{2^{2n + o(n)}}.$$
Taking the logarithm of both sides, we see that 
\begin{eqnarray*}
k \log n & > & \left(\frac{n-1}{2}\right) 
\log(n-1) - (2n + o(n)) \log 2 \\
& = & \frac{n}{2} \log(n-1) -\frac{1}{2} \log(n-1) - 2n \log 2 -
o(n) \\
        &  = & \frac{n}{2} \log n - o(n \log n) \\
        & = & (\frac{n}{2} -
          o(n)) \log n \end{eqnarray*} 
so $k > n/2 - o(n)$ as claimed. \endproof

It would be very interesting to find an upper bound for the SCR of any
polytope in $C_n$ that improves the $\bigo(n^2 \,\textup{log} \,n)$
upper bound on Chv{\'a}tal rank in \cite{EisenbrandSchulz}. Our
experiments in dimension up to $7$ suggest that there might be a
uniform upper bound for the SCR of any polytope in $C_n$ of order
$\bigo(n)$. Facet normals of $0/1$ $n$-polytopes with large
coefficients (matching the Alon-Vu bound) for $n \leq 10$ can be found
in the Polymake database at {\em
  http://www.math.tu-berlin.de/polymake/}. We have confirmed that for
$n \leq 7$, these facet normals appear in two rounds of IBN applied to
the normals of the standard relaxation of a $0/1$-polytope in $C_n$
used in the proof of Theorem~\ref{thm:0/1SCRlowerbound}. For instance,
when $n=7$, the Polymake database shows that $(9, 7, 5, 3, 2, 1, 1)$
is a possible facet normal. This vector lies in the minimal Hilbert
basis of the basis cone spanned by the vectors:
$$(3,2,2,2,1,0,0), (3,2,2,0,0,0,0), (3,2,1,0,0,0,0),(2,2,1,1,1,1,1),$$
$$(2,2,1,1,1,1,0), (2,2,1,1,1,0,1), (2,2,1,1,0,0,0)$$
which are all found in the first round of IBN applied to $\{\pm
1\}^7$. 

The fractional stable set polytope $\Frac(G)$ of a graph $G=([n],E)$
examined in Section~\ref{sec:stable set problem} lies in the unit cube
$C_n$. We will now derive a lower bound depending on $n$, for
$\scr(\Frac(G))$ as $G$ varies over all graphs with $n$ vertices. This
result contrasts the many examples of normals shown in
Section~\ref{sec:stable set problem} for which SCR is at most two. We
rely on a construction found in \cite{Liptak-Lovasz} for producing
facet normals of $\stab(G)$ with large coefficients.

\begin{definition} The {\bf product graph} of $G_1 = (V_1,E_1)$ and
  $G_2=(V_2,E_2)$ is the graph $G = (V,E)$ where $V = V_1 \cup V_2$
  and $E = E_1 \cup E_2 \cup \{uv \,:\, u \in V_1, v \in V_2 \}$.
\end{definition}

\begin{lemma} \label{lem:prod facet} Suppose $G_i=(V_i,E_i)$, $i=1,2$
  are graphs such that the inequality $\sum_{u \in V_i} a_i(u) x_u \leq b_i$ 
  defines a facet of $\stab(G_i)$ with $\aaa_i := (a_i(u))$ primitive
  for $i=1,2$. The inequality
$$ b_2 (\sum_{u \in V_1} a_1(u) x_u ) + b_1 (\sum_{u \in V_2} a_2(u)
x_u ) \leq b_1b_2$$ is facet-defining for $\stab(G)$. If $(b_1,b_2) =
1$ then the facet normal shown above is primitive.
\end{lemma}

\begin{observation} \label{obs:scr lower bound for frac} By
  Observation~\ref{obs:HAnorm}, for $G = ([n],E)$, if $\vv \in
  \A_G^{(k)}$, then $|| \vv ||_\infty \leq n^k$.  Therefore, if
  $\stab(G)$ has a primitive facet normal $\aaa \in \ZZ^n$, then
  $\scr(\Frac(G)) \geq \lfloor \textup{log}_n || \aaa ||_\infty
  \rfloor$.
\end{observation}

\begin{theorem} There is no constant $t$ such that $\scr(\Frac(G))
  \leq t$ for all graphs $G$.  
\end{theorem}

\begin{proof}
  Let $k_1, k_2, \ldots, k_p$ be the first $p$ prime numbers and
  consider the odd cycles $C_{2k_i+1}$ for $i=1,\ldots, p$. In each
  case, the odd hole inequality $\sum_{j=1}^{2k_i+1} x_j \leq k_i$ is
  facet-defining for $\stab(C_{2k_i+1})$. Let $G_p$ be the product
  graph $C_5 \times C_7 \times C_{11} \times \cdots C_{2k_p+1}$ which
  has $2(\sum_{i=1}^{p}k_i) + p$ vertices. By Lemma~\ref{lem:prod
    facet}, there is a primitive facet normal of $\stab(G_p)$ with
  infinity norm $\prod_{i=1}^{p} k_i$.
  
  The sum of the first $p$ prime numbers, $\sum_{i=1}^{p}k_i$ is
  approximately $\frac{1}{2}p^2 \textup{ln} p$ \cite{BachShallit},and
  hence the number of vertices of $G_p$ is approximately $p^2
  \,\textup{ln} \,p + p < p^{2+\epsilon}$. On the other hand,
  $\prod_{i=1}^{p} k_i$ is asymptotically $e^{(1+o(1)) p \ln p} 
  > e^{p \ln p}$. Therefore, $\scr(\Frac(G_p)) \geq
  \textup{log}_{p^{2+\epsilon}} e^{p \ln p} = \frac{p}{2+\epsilon}$,
  asymptotically.
 \end{proof}

\begin{problem} Is it true that for $G=([n],E)$, $\scr(\Frac(G)) \leq
  n$? More generally, is there an upper bound of order $\bigo(n)$ for
  the SCR of any polytope in the unit cube $C_n$?
\end{problem}

If the answer to the above problem is yes, then SCR would become
comparable to the number of steps needed by the modern {\em lift and
  project} methods for finding the integer hull of a polytope in $C_n$
such as those in \cite{BCC}, \cite{LovaszSchrijver91} and
\cite{SheraliAdams90}, since these methods take at most $n$
steps. There are a few different observations that support a positive
answer. For instance, it was shown in \cite{StephenTuncel} that the
semidefinite operator $N_+$ in \cite{LovaszSchrijver91} takes $\lfloor
n/2\rfloor$ iterations to produce $\stab(G)$ from $\Frac(G)$ when $G$
is the line graph of $K_n$ with $n$ odd. By Corollary~\ref{cor:line
  graphs}, $\scr(\Frac(G)) \leq 3$ for any line graph $G$. For the
operator $N$, it was shown in \cite{LovaszSchrijver91} that
$N(\Frac(G)) = Q_G^{(1)}$, the first Chv{\'a}tal closure of
$\Frac(G)$. Comparing with the small Chv{\'a}tal closure, we get that
$N^0(\Frac(G)) = \Frac(G) = S_G^{(0)}$ and $N^1(\Frac(G) = Q_G^{(1)}
\supseteq S_G^{(1)}$. If this pattern continues and we get
$N^k(\Frac(G)) \supseteq S_G^{(k)}$ for all $k \geq 2$, then indeed,
$\scr(\Frac(G))$ would be at most $n$ when $G$ has $n$ vertices.

%
%
\noindent{\bf Acknowledgments}. We thank Sasha Barvinok, 
Ravi Kannan and Les Trotter for helpful inputs to this paper.


\bibliographystyle{spmpsci}      


\end{document}